\documentclass[12pt]{amsart}

\usepackage[top=1in,left=1in,bottom=1in,right=1in]{geometry}
\usepackage{amsmath}
\usepackage{amssymb}
\usepackage{amsthm}
\usepackage{tikz}
\usetikzlibrary{arrows}

\usepackage{wrapfig}
\usepackage{graphicx}
\usepackage{epstopdf}
\usepackage{tikz-cd}

\usepackage{url}
\usepackage{scrextend}

\usepackage{colordvi}
\usepackage{color}

\usepackage{verbatim}

\usepackage{amsmath}
\usepackage{amsfonts}
\usepackage{colordvi}
\usepackage{color}
\usepackage{amsthm}
\usepackage{url}
\usepackage{graphicx}
\usepackage{epstopdf}
\usepackage{amssymb}

\newtheorem{innercustomthm}{Theorem}
\newenvironment{restatetheorem}[1]
  {\renewcommand\theinnercustomthm{#1}\innercustomthm}
  {\endinnercustomthm}
  
\newtheorem{theorem}{Theorem}[section]
\newtheorem{lemma}[theorem]{Lemma}

\newtheorem{proposition}[theorem]{Proposition}

\theoremstyle{definition}
\newtheorem{definition}[theorem]{Definition}
\theoremstyle{remark}
\newtheorem{remark}[theorem]{Remark}
\newtheorem{example}[theorem]{Example}

\DeclareMathOperator{\cone}{cone}

\DeclareMathOperator{\link}{lk}

\DeclareMathOperator{\del}{del}

\newcommand{\bergman}[1]{\underline{\Delta}_{#1}}
\newcommand{\indep}[1]{\mathcal I (#1)}
\newcommand{\augbergman}[1]{\Delta_{#1}}

\newcommand{\od}{:=}

\begin{document}
\title{Decompositions of Augmented Bergman Complexes}

\author{R. Amzi Jeffs}

\begin{abstract}
We study the augmented Bergman complex of a closure operator on a finite set, which interpolates between the order complex of proper flats and the independence complex of the operator. In 2020, Braden, Huh, Matherne, Proudfoot, and Wang showed that augmented Bergman complexes of matroids are always gallery-connected, and recently Bullock, Kelley, Reiner, Ren, Shemy, Shen, Sun, Tao, and Zhang strengthened ``gallery-connected" to ``shellable" by providing two classes of shelling orders: ``flag-to-basis" shellings and ``basis-to-flag" shellings. 

We show that augmented Bergman complexes of matroids are vertex decomposable, a stronger property than shellable. We also prove that the augmented Bergman complex of any closure operator is shellable if and only if lattice of flats (that is, its non-augmented Bergman complex) is shellable. As a consequence, an augmented Bergman complex is shellable if and only if it admits a flag-to-basis shelling. Perhaps surprisingly, the same does not hold for basis-to-flag shellings: we describe a closure operator whose augmented Bergman complex is shellable, but has no shelling order with bases appearing first. 
\end{abstract}

\thanks{Department of Mathematics, Carnegie Mellon University. Jeffs' work is supported by the National Science Foundation through Award No. 2103206.}
\date{\today}
\maketitle

\section{Introduction}

A \emph{closure operator} on a finite set $E$ is a function $f:2^E\to 2^E$ satisfying the following axioms for every $A,B\subseteq E$: \begin{itemize}
\item[C1.] $A\subseteq f(A)$, 
\item[C2.]  $A\subseteq B$ implies $f(A) \subseteq f(B)$, and
\item[C3.]  $f(f(A)) = f(A)$.
\end{itemize}
Each closure operator is determined by its lattice of \emph{flats} (sometimes called \emph{closed sets}), which is $\mathcal F(f) \od \{F\subseteq E\mid f(F) = F\}$. In particular, the meet of two flats is their intersection, and for any $A\subseteq E$, $f(A)$ is equal to the intersection of all flats containing $A$. 

To study the structure of a closure operator, one may associate it to three simplicial complexes: its Bergman complex $\bergman{f}$, its independence complex $\indep{f}$, and its augmented Bergman complex $\augbergman{f}$. The Bergman complex is the order complex of the lattice $\mathcal F(f)$, with the cone vertices $f(\emptyset)$ and $E$ removed, and the independence complex records subsets of $E$ whose closure gets smaller when any element is deleted. The augmented Bergman complex ``interpolates" between these two complexes, and in particular contains each as a full-dimensional induced subcomplex. Formally, these objects are defined as follows:
\begin{itemize}
\item The \emph{Bergman complex} $\bergman{f}$ has a vertex $x_F$ for every proper nonempty flat $F$ of $f$, and faces of the form $\{x_{F_1}, x_{F_2}, \ldots, x_{F_\ell}\}$ where $f(\emptyset)\subsetneq F_1\subsetneq F_2\subsetneq\cdots \subsetneq F_\ell\subsetneq E$ is a \emph{flag}, i.e. a chain of flats of $f$. We will sometimes refer to the \emph{cone over the Bergman complex}, denoted $\cone(\bergman{f})$, in which we allow $F_1 = f(\emptyset)$, or equivalently include the cone vertex $x_{f(\emptyset)}$. 
\item The \emph{independence complex} $\indep{f}$ has a vertex $y_i$ for every $i\in E$, and its faces are of the form $\{y_i\mid i\in I\}$ where $I\subseteq E$ is an \emph{independent set} of $f$, which means that $f(I\setminus \{i\}) \subsetneq f(I)$ for every $i\in I$. Independent sets with $f(I) = E$ are called \emph{bases}. 
\item The \emph{augmented Bergman complex} $\augbergman{f}$ has vertex set \[ \{y_i\mid i\in E\}\sqcup \{x_F\mid F\text{ is a proper flat of $f$}\}\] and its faces are \[
\{y_i\mid i\in I\} \sqcup \{x_{F_1}, \ldots, x_{F_\ell}\}
\]
where $I$ is an independent set of $f$, and $f(I)\subseteq F_1\subsetneq F_2\subsetneq \cdots \subsetneq F_\ell\subsetneq E$. 
\end{itemize}

Observe that $\cone(\bergman{f})$ is the induced subcomplex of $\augbergman{f}$ on the $x_F$ vertices, and $\indep{f}$ is the induced subcomplex of $\augbergman{f}$ on the $y_i$ vertices. 

An important special case is when $f$ is the closure operator of a matroid $M$. In this case $\indep{f}$ is the usual independence complex of $M$, and $\bergman{f}$ is the usual Bergman complex of $M$, i.e. the order complex of the proper part of its lattice of flats. We will use $\augbergman{M}$ to denote the augmented Bergman complex in this case. Augmented Bergman complexes of matroids played an important role in the recent foundational work of Braden, Huh, Matherne, Proudfoot, and Wang \cite{semismall, singularhodge}, in which it was shown that they are connected in codimension-1 (or ``gallery connected"). Recently, Bullock, Kelley, Reiner, Ren, Shemy, Shen, Sun, Tao, and Zhang \cite{augbergman} strengthened this result by showing that the augmented Bergman complex of a matroid is always shellable. In particular, they showed that the augmented Bergman complex of a matroid admits two classes of shelling orders: \emph{flag-to-basis} shellings in which maximal flags appear first and bases appear last, and \emph{basis-to-flag} shellings in which the reverse occurs. 

We will show that the augmented Bergman complex of a matroid is vertex decomposable, a stronger property than shellability. For formal definitions of shellability and vertex decomposability, see Section \ref{sec:background}. 

\begin{theorem}\label{thm:matroidvd}
The augmented Bergman complex of a matroid is vertex decomposable. 
\end{theorem}

The proof of Theorem \ref{thm:matroidvd} proceeds by deleting vertices corresponding to flats according to a linear extension of the lattice of flats. One side effect of this strategy is that we obtain a variety of induced subcomplexes of the augmented Bergman complex that are vertex decomposable. Our proof depends on the matroid structure of $M$, and the fact that the independence complex and Bergman complex of a matroid are both vertex decomposable. See Proposition \ref{prop:upperset} for details. 

We also show, for any closure operator $f$, that shellability of $\augbergman{f}$ is completely determined by shellability of $\bergman{f}$.  
In particular, whenever $\bergman{f}$ is shellable we construct a shelling order of $\augbergman{f}$ with maximal flags appearing first and bases last. 
Our shelling orders are a very slight generalization of the flag-to-basis shellings constructed by \cite{augbergman} in the matroidal case, see Remark \ref{rem:flagtobasistechnicality}.
Below, $f/F$ is the \emph{contraction} of $f$ by a proper flat $F$, defined formally in Section \ref{sec:background}. 

\begin{theorem}\label{thm:shellability}
Let $f$ be a closure operator on a finite set $E$. The following are equivalent:\begin{itemize}
\item[(i)] The Bergman complex $\bergman{f}$ is shellable,
\item[(ii)] For every proper flat $F$ of $f$, the Bergman complex $\bergman{f/F}$ is shellable,
\item[(iii)] The augmented Bergman complex $\augbergman{f}$ admits a shelling order with maximal flags appearing first and bases appearing last, and
\item[(iv)] The augmented Bergman complex $\augbergman{f}$ is shellable. 
\end{itemize}
\end{theorem}

Given the equivalence of (iii) and (iv) above, and the results of \cite{augbergman}, one might also expect that $\augbergman{f}$ is shellable if and only if it admits a basis-to-flag shelling. However, the following example refutes this by providing a closure operator whose bases do not generate a shellable complex. In fact, the bases of this operator are exactly the maximal independent sets, so~$\augbergman{f}$ is shellable while $\indep{f}$ is not. 

\begin{example}\label{ex:intro}
Let $E = \{1,2,3,4,5\}$ and let $f$ be the closure operator whose proper flats are the empty set, all singleton sets, and the pairs $\{1,2\}, \{1,3\}, \{2,3\}, \{3,4\}, \{3,5\}$, and $\{4,5\}$. One may check that the independence complex of $f$ consists of a complete graph on $E$, plus the triangles $\{1,2,3\}$ and $\{3,4,5\}$, and the maximal independent sets are exactly the bases. Note that $\indep{f}$ is not shellable because the two triangles share only a single vertex. 

Figure \ref{fig:no-converse} shows the augmented Bergman complex of $f$, in three layers. The bottom layer is $\cone(\bergman{f})$, the top layer is $\indep{f}$, and the middle layer consists of the ``hybrid" faces in~$\augbergman{f}$, which contain vertices corresponding to both flats and ground set elements. One can form a shelling order of $\augbergman{f}$ by first shelling the bottom layer, then adding the facets from the middle layer which share an edge with the bottom, then adding the facets in the middle layer which share only one vertex with the bottom layer, and finally adding the facets from the top layer in any order. 
Our proof of Theorem \ref{thm:shellability} will generalize this approach to shelling the augmented Bergman complex. 
\begin{figure}
\[
\includegraphics{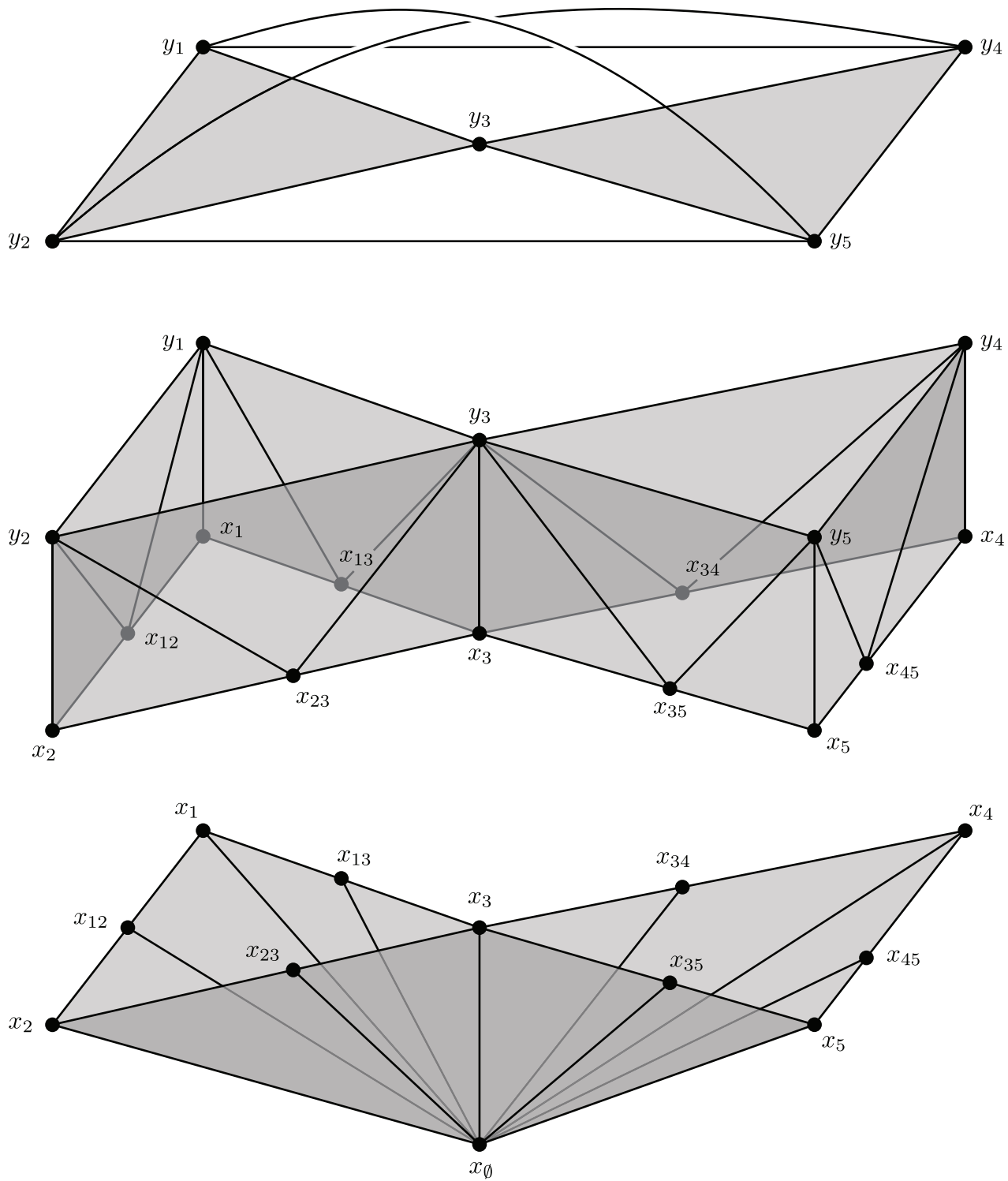}
\]
\caption{The augmented Bergman complex of the operator $f$ from Example \ref{ex:intro}, broken into three layers.}\label{fig:no-converse}
\end{figure}
\end{example}

\section{Background}\label{sec:background}

We first recall some general definitions and notation regarding simplicial complexes. A simplicial complex $\Delta$ is \emph{shellable} its facets can be ordered $\sigma_1, \sigma_2, \ldots, \sigma_k$ so that for $2\le i \le k$ the simplicial complex \[
\langle \sigma_1, \sigma_2, \ldots, \sigma_{i-1}\rangle \cap \langle \sigma_i\rangle
\]
is pure of dimension $\dim(\sigma_i)-1$. Above, $\langle \sigma_1, \sigma_2, \ldots, \sigma_{i-1}\rangle$ denotes the simplicial complex generated by a collection of faces. 
The \emph{deletion} of a face $\sigma$ in $\Delta$ is the simplicial complex\[
\del_\Delta(\sigma) \od \{\tau \setminus \sigma \mid \tau \in \Delta\}.
\]
The \emph{link} of a face $\sigma$ in $\Delta$ is the simplicial complex\[
\link_\Delta(\sigma) \od \{\tau \setminus \sigma  \mid \sigma\subseteq \tau \in \Delta\}.
\]
If $\Delta$ and $\Gamma$ are simplicial complexes on disjoint vertex sets, the \emph{join} of $\Delta$ and $\Gamma$ is the simplicial complex \[
\Delta \ast \Gamma \od \{\sigma \cup \tau \mid \sigma \in \Delta\text{ and } \tau \in \Gamma\}. 
\]

\begin{definition}\label{def:vertexdecomposable}
A simplicial complex $\Delta$ is \emph{vertex decomposable} if $\Delta$ is a simplex (including the possibility $\Delta = \{\emptyset\}$) or, inductively, if there is a vertex $v$ of $\Delta$ so that \begin{itemize}
\item[(i)] $\del_\Delta(v)$ and $\link_\Delta(v)$ are both vertex decomposable, and
\item[(ii)] every facet of $\del_\Delta(v)$ is also a facet of $\Delta$.
\end{itemize}
A vertex satisfying these conditions is called a \emph{decomposing vertex}. A vertex that satisfies~(ii) is called a \emph{shedding vertex}. 
\end{definition}

Note that $v$ is a shedding vertex if and only if the following holds: for each facet $\sigma$ containing $v$, there is another facet $\tau$ with $\sigma\setminus \tau = \{v\}$.

Every vertex decomposable complex is shellable. Indeed, given a decomposing vertex $v$ in~$\Delta$, one may form a shelling order 
\[
\sigma_1, \sigma_2, \ldots, \sigma_k, \tau_1\cup \{v\}, \tau_2\cup \{v\}, \ldots,  \tau_\ell \cup\{v\}
\]
of $\Delta$, where the $\sigma_i$ are a shelling of $\del_\Delta(v)$ and the $\tau_j$ are a shelling of $\link_\Delta(v)$. Furthermore, if $\Delta$ and $\Gamma$ are vertex decomposable, then so is $\Delta\ast \Gamma$, a fact that we will make use of later. 

Now let us establish some conventions and definitions that are specific to augmented Bergman complexes and closure operators. If $F$ is a flat of a closure operator $f$, the \emph{contraction} of $f$ by $F$ is the closure operator $f/F: 2^{E\setminus F} \to 2^{E\setminus F}$ with $(f/F)(A) \od f(A\cup F) \setminus F$. The \emph{restriction} of $f$ to $F$ is the closure operator $f|_F: 2^{F} \to 2^{F}$ defined by $(f|_F)(A) \od f(A)$. Observe that the flats of $f/F$ are exactly the flats of $f$ that contain $F$, but with the elements of $F$ removed. Moreover, $\indep{f|_F}$ consists of the faces of $\indep{f}$ that are contained in $F$. Lastly, note that when $f$ is the closure operator of a matroid $M$, $f/F$ is the closure operator arising from the contraction of $M$ by $F$ (denoted $M/F$), and $f|_F$ is the closure operator arising from the restriction of $M$ to $F$ (denoted $M|_F$). 

Following \cite{augbergman}, we will denote faces of the augmented Bergman complex by pairs $\phi = (I, F_\bullet)$ where $I$ is an independent set and $F_\bullet$ denotes a flag $F_1\subsetneq F_2\subsetneq\cdots \subsetneq F_\ell$ of proper flats that is compatible with $I$ in the sense that $I\subseteq F_1$. We sometimes abuse notation and regard $F_\bullet$ as a flag of proper flats of $f/F_1$. Observe that if $\phi$ is a facet, then $f(I) = F_1$, or $F_\bullet$ is empty and $f(I) = E$. 

We need one last observation regarding augmented Bergman complexes. The following lemma is straightforward, but plays a crucial role in our later proofs. 

\begin{lemma}\label{lem:linkofflat}
Let $F$ be a proper flat of a closure operator $f$. Then \[
\link_{\augbergman{f}}(x_F) \cong \augbergman{f|_F} \ast \bergman{f/F}. 
\]
In words, the link of $x_F$ in the augmented Bergman complex of $f$ is the join of the augmented Bergman complex of $f|_F$ with the Bergman complex of $f/F$. 
\end{lemma}
\begin{proof}
The faces $(I, F_\bullet)$ of $\link_{\bergman{f}}(x_F)$ are exactly those for which $I\subseteq F$ and $F$ can be inserted into the chain $F_\bullet$. We may write these faces uniquely as $(I, F_\bullet') \sqcup F_\bullet''$ where $F_\bullet'$ consists of the flats in $F_\bullet$ that are properly contained in $F$, and $F_\bullet''$ consists of the flats in $F_\bullet$ that properly contain $F$. Note that the pairs $(I,F_\bullet')$ are exactly the faces of the augmented Bergman complex of $f|_F$, while the chains $F_\bullet''$ are exactly the faces of $\bergman{f/F}$. This proves the result. \end{proof}

\section{Vertex Decomposability}\label{sec:vd}

We are now ready to prove Theorem \ref{thm:matroidvd}. In fact, we will prove a somewhat stronger result, arguing that the augmented Bergman complex of a matroid has a variety of induced subcomplexes that are vertex decomposable, each with concrete choices of decomposing vertex. Below, an \emph{upper-set} of proper flats $\mathcal L$ is a collection of proper flats so that if $F\in \mathcal L$ and $F'$ is a proper flat containing $F$, then $F'\in \mathcal L$. 

Note that Theorem \ref{thm:matroidvd} follows from Proposition \ref{prop:upperset} by choosing $\mathcal L = \mathcal F(M)\setminus \{E\}$. 

\begin{proposition}\label{prop:upperset}
Let $M$ be a matroid on ground set $E$. Let~$\mathcal L$ be an upper-set of proper flats, and let $\augbergman{M}(\mathcal L)$ be the subcomplex of $\augbergman{M}$ induced on the vertex set $\{y_i\mid i\in E\}\sqcup \{x_F\mid F \in \mathcal L\}$. Then $\augbergman{M}(\mathcal L)$ is vertex decomposable. Moreover, if $F_0$ is a minimal element of $\mathcal L$, then $x_{F_0}$ is a decomposing vertex of $\augbergman{M}(\mathcal L)$.
\end{proposition}
\begin{proof}

We work by induction on the size of $\mathcal L$. When $\mathcal L$ is empty, $\Delta_M(\mathcal L)$ is just $\indep{M}$. The independence complex of a matroid is always vertex decomposable (Provan and Billera \cite{provanbillera} observed that any vertex is a decomposing vertex) and so the result follows in this case. For the inductive step, suppose that $\mathcal L$ is nonempty and let $F_0$ be a minimal element of $\mathcal L$. The deletion $\del_{\augbergman{M}(\mathcal L)}(x_{F_0})$ is equal to $\augbergman{M}(\mathcal L \setminus \{x_{F_0}\})$, which is vertex decomposable by inductive hypothesis. 

By Lemma \ref{lem:linkofflat}, the link of $x_{F_0}$ in $\augbergman{M}$ is $\augbergman{M|_{F_0}}\ast \bergman{M/F_0}$. In $\augbergman{M}(\mathcal L)$ we have deleted all vertices corresponding to flats contained in $F_0$, but no vertices corresponding to flats containing $F_0$. Thus when considering the link of $x_{F_0}$ in $\augbergman{M}(\mathcal L)$ the first term in the join above becomes $\indep{M|_{F_0}}$ while the second term is unaltered, and we have \[
\link_{\augbergman{M}(\mathcal L)}(x_{F_0})\cong \indep{M|_{F_0}}\ast \bergman{M/F_0}. 
\]
As noted above, the independence complex of a matroid is vertex decomposable. Moreover, the lattice of flats of a matroid admits a CL-labeling (see \cite[Section 7.6]{bjornerchapter}), which implies that its order complex (i.e. the Bergman complex of the matroid) is vertex decomposable (see \cite[Theorem 11.6]{bjornerwachs2}). Thus both terms in the join above are vertex decomposable, and we conclude that $\link_{\augbergman{M}(\mathcal L)}(x_{F_0})$ is vertex decomposable. 

 It remains to argue that $x_{F_0}$ is a shedding vertex. Let $(I, F_\bullet)$ be a facet of $\Delta_M(\mathcal L)$ that contains $x_{F_0}$. Then $I\subseteq F_0$, and $F_0$ is equal to the first element of the chain $F_\bullet$. If $F_0$ is the only flat in $F_\bullet$ then---because $M$ is a matroid---we may choose $a\in E\setminus I$ so that $I\cup \{a\}$ is an independent set. Then the facet $(I\cup \{a\}, \emptyset)$ contains all vertices of $(I, F_\bullet)$ except for~$x_{F_0}$. If $F_\bullet$ is a chain $F_0 = F_1\subsetneq F_2 \subsetneq \cdots \subsetneq F_\ell$, then---again, because $M$ is a matroid---we may choose $a\in F_2\setminus I$ so that $I\cup \{a\}$ is an independent set. Define $F_\bullet'$ to be the flag $F_2\subsetneq \cdots \subsetneq F_\ell$, and observe that the facet corresponding to $(I\cup \{a\}, F_\bullet')$ contains all vertices of $(I, F_\bullet)$ except for $x_{F_0}$, proving the result.  
\end{proof}

\section{Shellability}\label{sec:shellability}

We now proceed to our proof of Theorem \ref{thm:shellability}. The main substance of the argument lies in proving that (ii) implies (iii). We apply techniques similar to \cite{augbergman} to construct the desired shelling order, but the details of our proof are somewhat different since we are not working with the closure operator of a matroid. 

\begin{restatetheorem}{\ref{thm:shellability}}
Let $f$ be a closure operator on a finite set $E$. The following are equivalent:\begin{itemize}
\item[(i)] The Bergman complex $\bergman{f}$ is shellable,
\item[(ii)] For every proper flat $F$ of $f$, the Bergman complex $\bergman{f/F}$ is shellable,
\item[(iii)] The augmented Bergman complex $\augbergman{f}$ admits a shelling order with maximal flags appearing first and bases appearing last, and
\item[(iv)] The augmented Bergman complex $\augbergman{f}$ is shellable. 
\end{itemize}
\end{restatetheorem}
\begin{proof}
Clearly (iii) implies (iv).
Furthermore, (iv) implies (i) since the link of $x_{f(\emptyset)}$ in $\augbergman{f}$ is $\bergman{f}$, and shellability is inherited by links. 
To see that (i) implies (ii), first note that the link of $x_F$ in $\bergman{f}$ is $\bergman{f|_F} \ast \bergman{f/F}$, and a join of two complexes is shellable if and only if both complexes are shellable (see \cite[Remark 10.22]{bjornerwachs2}). 
Thus shellability of $\bergman{f}$ implies shellability of $\bergman{f|_F}\ast \bergman{f/F}$ for every proper flat, which in turn implies shellability of $\bergman{f/F}$ as desired.
It remains to argue that (ii) implies (iii). 

Suppose that (ii) holds, and fix a shelling order for every Bergman complex $\bergman{f/F}$ where~$F$ is a proper flat.
Let $<$ be a linear extension of the independence complex $\indep{f}$---that is, $<$ is a total order on independent sets so that $I<I'$ whenever $I\subseteq I'$. 
Then define a total order~$\prec$ on the facets of $\augbergman{f}$ as follows.
If $\phi = (I, F_\bullet)$ and $\phi' = (I', F_\bullet')$ both have nonempty flags, then $\phi \prec \phi'$ whenever \begin{itemize}
\item[(a)] $I<I'$, or
\item[(b)] $I = I'$ and $F_\bullet\setminus \{x_{f(I)}\}$ precedes $F_\bullet'\setminus \{x_{f(I)}\}$ in our fixed shelling of $\bergman{f/f(I)}$.
\end{itemize}
Extend $\prec$ to facets with empty flags by placing them after the facets with nonempty flags, in any order. Observe that $\prec$ begins with with facets of the form $(\emptyset, F_\bullet)$, i.e. maximal flags, and ends with facets of the form $(I, \emptyset)$, i.e. bases. 

We must argue that $\prec$ is a shelling order. This amounts to showing that for every $\phi$ and~$\phi'$ with $\phi \prec \phi'$ we can construct $\phi''$ so that $\phi'' \prec \phi'$, $\phi\cap \phi'\subseteq \phi''$, and $|\phi'\setminus \phi''| = 1$. We consider two cases. 

\emph{Case 1:} Suppose that $I\neq I'$.
Our ordering guarantees that $I'$ is not contained in $I$, so we may choose $a\in I'\setminus I$. 
Define $I'' = I'\setminus \{a\}$, and note that $f(I'')$ is a flat properly contained in~$f(I')$. 
Hence we may choose a flag $F_\bullet''$ which is maximal among flags whose minimal element is $f(I'')$ and which contain $F_\bullet'$. 
The facet $\phi'' = (I'', F_\bullet'')$ precedes $\phi'$ since $I'' < I'$, and it contains every vertex of $\phi'$ except $y_a$. 
Since $y_a$ is not a vertex of $\phi$, we see that $\phi''$ contains $\phi\cap \phi'$.

\emph{Case 2:} Suppose that $I = I'$.
Then $F_\bullet\setminus \{x_{f(I)}\}$ precedes $F_\bullet'\setminus  \{x_{f(I)}\}$ in our fixed shelling of $\augbergman{f/f(I)}$, and we may choose a facet of $\augbergman{f/f(I)}$ that precedes the latter facet, contains all but one of its vertices, and contains the intersection of these two facets.
Let $F_\bullet''$ be the result of adding $x_{f(I)}$ to this facet, and note that $\phi'' = (I'', F_\bullet'')$ is a facet of $\augbergman{f}$, where $I'' = I$. 
By construction, $\phi''$ precedes $\phi'$, contains $\phi\cap \phi'$, and contains all but one vertex of $\phi'$.
We conclude that $\prec$ is a shelling order, and the theorem follows. 
\end{proof}

\begin{remark}\label{rem:flagtobasistechnicality}
Conditions (a) and (b) above are analogous to the conditions in the definition of a flag-to-basis shelling given in \cite[Definition 3.1]{augbergman}. 
Our conditions are very slightly more general, in that we allow for linear extensions of the independence complex that are not necessarily monotone in the size of independent sets.
However, if one restricts to linear extensions that are monotone in the size of independent sets, then our shelling order is exactly a flag-to-basis shelling. 
\end{remark}

\begin{remark}\label{rem:restrictionsets}
In any shelling order for a simplicial complex, each facet $F$ can be associated to its \emph{restriction set}, denoted $\mathcal R(F)$, which is the unique minimal face contained in $F$ but no previous facet. 
The shelling order $\prec$ that we constructed in the proof of Theorem \ref{thm:shellability} has the advantage that its restriction sets can be described succinctly. 
One may verify that for a facet $\phi = (I, F_\bullet)$ with $F_\bullet$ nonempty, we have $\mathcal R(\phi) = I\sqcup \mathcal R(F_\bullet\setminus \{x_{F_1}\})$, where $\mathcal R(F_\bullet\setminus \{x_{F_1}\})$ is the restriction set of $F_\bullet\setminus \{x_{F_1}\}$ in the fixed shelling of $\bergman{f/F_1}$ used to define $\prec$. Furthermore, if $F_\bullet$ is empty, then $\mathcal R(\phi)$ is simply $I$. 

For pure shellable complexes, restriction sets can be used to compute the \emph{$h$-vector}, an important numerical invariant of the complex. 
In this context it turns out that $h_i$, the $i$-th entry of the $h$-vector, is exactly the number of facets whose restriction set has size~$i$ (see Section 2 of Chapter III in \cite{combcommalg}). 
Often it is convenient to record the $h$-vector of a $d$\nobreakdash-dimensional complex $\Delta$ by the \emph{$h$-polynomial} $h(\Delta, t) \od \sum_{i=0}^d h_i t^i$. 
Our observation about the restriction sets of the shelling order above yields the following formula for the $h$-polynomial of a shellable augmented Bergman complex, provided that it is additionally pure (which happens, for example, when $f$ is the closure operator of a matroid):
\[
h(\augbergman{f}, t) = \sum_{I\in \indep{f}} t^{|I|} h(\bergman{f/f(I)}, t). 
\]
Above, $h(\bergman{f/f(I)}, t) = 1$ whenever $I$ is a basis.

We conclude by noting that this formula can fail when $\augbergman{f}$ is not pure.
The closure operator from Example \ref{ex:intro} has an augmented Bergman complex that is shellable, but whose $h$-vector is equal to $(1,14,19,-2)$. 
However, the Bergman complexes of its contractions by proper flats are all pure and shellable, so the formula above would give nonnegative coefficients for its $h$-polynomial. 
Thus the formula above does not hold for this closure operator, despite the fact that its augmented Bergman complex is shellable. 
It may be possible to obtain similar formulas for nonpure shellable complexes using the $h$-triangle defined by Bj\"orner and Wachs \cite[Section 3]{bjornerwachs1}.
\end{remark}

\section*{Acknowledgements}
We thank Vic Reiner for useful discussions and feedback on this work. 
We are also grateful to Jos\'e Samper for suggesting that we investigate the structure of restriction sets, leading to our comments in Remark \ref{rem:restrictionsets}. 

\bibliographystyle{plain}
\bibliography{references}

\begin{thebibliography}{1}

\bibitem{bjornerwachs1}
Anders Bj\"orner and Michelle~L. Wachs.
\newblock Shellable nonpure complexes and posets. {I}.
\newblock {\em Transactions of the American Mathematical Society},
  348(4):1299--1327, 1996.

\bibitem{bjornerwachs2}
Anders Bj\"orner and Michelle~L. Wachs.
\newblock Shellable nonpure complexes and posets. {II}.
\newblock {\em Transactions of the American Mathematical Society},
  349(10):3945–3975, 1997.

\bibitem{bjornerchapter}
Anders Björner.
\newblock {\em Homology and Shellability of Matroids and Geometric Lattices},
  page 226–283.
\newblock Encyclopedia of Mathematics and its Applications. Cambridge
  University Press, 1992.

\bibitem{semismall}
Tom Braden, June Huh, Jacob~P. Matherne, Nicholas Proudfoot, and Botong Wang.
\newblock A semi-small decomposition of the chow ring of a matroid, 2020.
\newblock https://arxiv.org/abs/2002.03341.

\bibitem{singularhodge}
Tom Braden, June Huh, Jacob~P. Matherne, Nicholas Proudfoot, and Botong Wang.
\newblock Singular hodge theory for combinatorial geometries, 2020.
\newblock https://arxiv.org/abs/2010.06088.

\bibitem{augbergman}
Elisabeth Bullock, Aidan Kelley, Victor Reiner, Kevin Ren, Gahl Shemy, Dawei
  Shen, Brian Sun, and Zhichun~Joy Zhang.
\newblock Topology of {A}ugmented {B}ergman {C}omplexes.
\newblock {\em Electron. J. Combin.}, 29(1):Paper No. 1.31, 2022.

\bibitem{provanbillera}
J.~Scott Provan and Louis~J. Billera.
\newblock Decompositions of simplicial complexes related to diameters of convex
  polyhedra.
\newblock {\em Mathematics of Operations Research}, 5(4):576--594, 1980.

\bibitem{combcommalg}
Richard~P. Stanley.
\newblock {\em Combinatorics and commutative algebra}, volume~41 of {\em
  Progress in Mathematics}.
\newblock Birkh\"{a}user Boston, Inc., Boston, MA, second edition, 1996.

\end{thebibliography}

\end{document}